\documentclass[a4paper,10pt]{amsart}

\usepackage{amssymb,amsmath,amsfonts,amsthm,mathtools}
\usepackage{latexsym,mathrsfs,pb-diagram}
\usepackage[pdftex]{graphicx}
\usepackage[all]{xy}
\usepackage[hmargin=2.5cm,vmargin=3.5cm]{geometry}
\usepackage[plainpages=false,colorlinks,pdfpagelabels]{hyperref}
\hypersetup{
  urlcolor=black,
  citecolor=green,
  linkcolor=blue
}

\numberwithin{equation}{section}

\theoremstyle{definition}
\newtheorem{Def}{Definition}[section]

\theoremstyle{remark}

\theoremstyle{plain}
\newtheorem{Prop}[Def]{Proposition}
\newtheorem{Cor}[Def]{Corollary}
\newtheorem{Thm}[Def]{Theorem}

\newcommand{\dfn}{\mathrel{\dot{=}}}

\newcommand{\st}{ \ ; \ }
\newcommand{\rarr}{\rightarrow}
\newcommand{\sset}{\subset}

\newcommand{\Z}{\mathbb{Z}}
\newcommand{\N}{\mathbb{N}}
\newcommand{\R}{\mathbb{R}}

\newcommand{\C}{\mathbb{C}}

\newcommand{\TR}[5]{\begin{array}{c c c c c}
    {#1} & : & {#3} & \longrightarrow & {#5}\\
    & & {#2} & \longmapsto & {#4}
  \end{array}
}

\newcommand{\transp}[1]{\prescript{\mathrm{t}}{}{{#1}}}

\DeclareMathOperator{\Span}{\mathrm{span}}
\DeclareMathOperator{\ran}{\mathrm{ran}}

\newcommand{\del}{\partial}
\newcommand{\dd}{\mathrm{d}}
\renewcommand{\Re}{\mathsf{Re}}

\newcommand{\cinfty}{\mathscr{C}^\infty}

\newcommand{\MM}{\mathrm{M}}
\newcommand{\LL}{\mathrm{L}}
\newcommand{\VV}{\mathcal{V}}

\newcommand{\e}{\mathbf{e}}

\newcommand{\s}{\mathbb{S}}

\DeclareMathOperator{\SU}{\mathrm{SU}}
\DeclareMathOperator{\GL}{\mathrm{GL}}
\DeclareMathOperator{\su}{\mathfrak{su}}
\DeclareMathOperator{\gl}{\mathfrak{gl}}

\DeclareMathOperator{\End}{\mathrm{End}}
\DeclareMathOperator{\Aut}{\mathrm{Aut}}

\newcommand{\gr}[1]{\mathfrak{#1}}

\DeclareMathOperator{\ad}{\mathrm{ad}}
\newcommand{\vv}[1]{\mathrm{#1}}

\author{Gabriel Ara\'{u}jo}
\address{University of S{\~a}o Paulo, ICMC-USP, S{\~a}o Carlos, SP, Brazil}
\email{\texttt{gccsa@icmc.usp.br}}
\thanks{This work was partially supported by Conselho Nacional de Desenvolvimento Cient{\'i}fico e Tecnol{\'o}gico (CNPq, grant~140838/2012-0) and the S{\~a}o Paulo Research Foundation (FAPESP, grant~2018/12273-5).}

\title[Involutive structures on $\SU(2)$: examples]{Computing cohomology spaces of left-invariant involutive structures on $\SU(2)$: examples}

\begin{document}

\begin{abstract} In these notes we study left-invariant involutive structures on $\SU(2)$, the most na{\"i}ve non-commutative compact Lie group. We determine closedness of the range (in the smooth topology) of a single complex vector field spanning the standard CR structure of $\SU(2)$ and also compute the smooth cohomology spaces of a corank~$1$ structure. In our approach, it is fundamental to understand concretely the irreducible representations of the ambient Lie group and how left-invariant vector fields operate on their matrix coefficients (which we borrow from~\cite[Chapter~11]{rt_psdoas} and briefly recall below).

  Our purpose is solely to provide some easy applications of the theory developed in~\cite{araujo19}, as the results shown here can probably be obtained by more direct methods.
\end{abstract}

\maketitle


\section{General remarks and notation}

Let $N \in \N$. By $\gl(N, \C)$ we denote the (complex) vector space of all $N \times N$ matrices with complex entries. A matrix in $\gl(N, \C)$ acts on $\C^N$ via left multiplication, yielding a $\C$-linear endomorphism of $\C^N$. Every such endomorphism has this form; this establishes an isomorphism of complex vector spaces $\gl(N, \C) \cong \End(\C^N)$. Moreover, both these spaces carry natural algebra structures (matrix multiplication in $\gl(N, \C)$ and composition of linear operators in $\End(\C^N)$), hence Lie algebra structures (the commutators for the underlying algebras), for which the aforementioned linear isomorphism is also a homomorphism: $\gl(N, \C)$ and $\End(\C^N)$ are also isomorphic as Lie algebras.

Let $\langle \cdot, \cdot \rangle$ denote the standard Hermitian inner product on $\C^N$, that is 
\begin{align*}
  \langle z, w \rangle &\dfn \sum_{j = 1}^N z_j \bar{w}_j
\end{align*}
whenever $z = (z_1, \ldots, z_N), w = (w_1, \ldots, w_N) \in \C^N$. For $A \in \End(\C^N)$ we denote by $A^* \in \End(\C^N)$ its adjoint w.r.t.~$\langle \cdot, \cdot \rangle$ i.e.
\begin{align*}
  \langle A z, w \rangle &= \langle z, A^* w \rangle, \quad \forall z,w \in \C^N.
\end{align*}
In terms of the corresponding matrix we have
\begin{align*}
  A = (a_{jk})_{1 \leq j,k \leq N} &\Longrightarrow A^* = (\bar{a}_{kj})_{1 \leq j,k \leq N}
\end{align*}
that is, $A^*$ is simply the conjugate transpose $\transp{\bar{A}}$. Notice also that, as a simple computation shows, the usual Euclidean inner product  (the dot product) of $z, w \in \C^N \cong \R^{2N}$can be written as $\Re \langle z, w \rangle$.

We denote by $\GL(N, \C)$ the set of all invertible matrices in $\gl(N, \C)$. Via the isomorphism of algebras $\gl(N, \C) \cong \End(\C^N)$, it is identified with $\Aut(\C^N)$, the set of all complex linear automorphisms of $\C^N$. Both of them carry natural Lie group structures, turning the identification $\GL(N, \C) \cong \Aut(\C^N)$ into a Lie group isomorphism.

\section{The special unitary group}

Our object of interest is the \emph{special unitary group} of order $2$:
\begin{align*}
  \SU(2) &\dfn \{ x \in \GL(2, \C) \st x^* x = e, \ \det x = 1 \}
\end{align*}
which is clearly a closed (hence Lie) subgroup of $\GL(2, \C)$. Here, $e$ stands for the $2 \times 2$ identity matrix. Consider the map $\theta: \C^2 \rarr \gl(2, \C)$ defined by
\begin{align*}
  \theta(z_1,z_2) &\dfn \left(
    \begin{array}{c c}
      z_1 & -\bar{z}_2 \\
      z_2 & \bar{z}_1
    \end{array}
    \right)
\end{align*}
This map is clearly $\R$-linear and injective, hence $\ran \theta$ is a $\R$-linear $4$-dimensional subspace of $\gl(2, \C)$. It is easy to check that every $x \in \SU(2)$ is of the form $x = \theta(z)$ for some $z \in \C^2$ satisfying $|z| = 1$ i.e.~$\SU(2)$ is the image under $\theta$ of the $3$-sphere
\begin{align*}
  \s^3 &\dfn \{ z \in \C^2 \st |z| = 1 \}
\end{align*}
hence diffeomorphic to it. Since $\theta$ is $\R$-linear and injective it maps tangent vectors to $\s^3$ to tangent vectors to $\SU(2)$. It is easy to see that
\begin{align*}
  T_{(1,0)} \s^3 &= i\R \times \C
\end{align*}
(use curves in $\s^3$ passing through $(1,0)$) and since $\theta(1,0) = e$ we have
\begin{align*}
  T_e \SU(2) = \theta (i\R \times \C) = \left\{
  \left(
    \begin{array}{c c}
      it & -\bar{z}_2 \\
      z_2 & -it
    \end{array}
    \right) \st t \in \R, \ z_2 \in \C \right\}
\end{align*}
which is easily seen to be precisely
\begin{align*}
  \su(2) &\dfn \{ \vv{X} \in \gl(2, \C) \st \vv{X} + \vv{X}^* = 0, \ \mathrm{tr} \vv{X} = 0 \}.
\end{align*}
One can prove that $\su(2)$ is a Lie subalgebra of $\gl(2, \C)$ and that the map ``evalutation at the identity''
\begin{align*}
  \text{$\vv{X}$ a real left-invariant vector field on $\SU(2)$} &\longmapsto \vv{X}(e) \in T_e \SU(2) \cong \su(2)
\end{align*}
is a Lie algebra isomorphism. For that reason, one regards $\su(2)$ as the Lie algebra of $\SU(2)$: a matrix $\vv{X} \in \su(2)$ acts on $f \in C^\infty(\SU(2))$ via the formula
\begin{align*}
  (\vv{X} f)(x) &\dfn \left. \frac{\dd}{\dd t} \right|_{t = 0} f \left(x \cdot \e^{t \vv{X}} \right), \quad x \in \SU(2),
\end{align*}
where
\begin{align*}
  \e^{t \vv{X}} &= \sum_{k = 0}^\infty \frac{t^k \vv{X}^k}{k!}
\end{align*}
is proven to belong to $\SU(2)$ for every $t \in \R$ provided $\vv{X} \in \su(2)$.

We are interested in studying left-invariant involutive structures on $\SU(2)$. As we have seen, these are Lie subalgebras of $\C \su(2)$ (which, by the way, is precisely $\gr{sl}(2, \C)$ -- but we shall not use this fact). For this purpose, next we borrow some computations and notations from~\cite[Chapter~11]{rt_psdoas}.
\subsection{Choosing convenient frames and metrics}

We define
\begin{align*}
  \vv{Y}_1 &\dfn \frac{1}{2} \left( \begin{array}{c c} 0 & i \\ i & 0 \end{array} \right) = \theta(0, i/2), \\
  \vv{Y}_2 &\dfn \frac{1}{2} \left( \begin{array}{c c} 0 & -1 \\ 1 & 0 \end{array} \right) = \theta(0,1/2), \\
  \vv{Y}_3 &\dfn \frac{1}{2} \left( \begin{array}{c c} i & 0 \\ 0 & -i \end{array} \right) = \theta(i/2,0)
\end{align*}
which clearly forms a basis for $\su(2)$. For more convenient computations in the forthcoming sections we shall further define
\begin{align*}
  \del_+ &\dfn i \vv{Y}_1 - \vv{Y}_2 \\
  \del_- &\dfn i \vv{Y}_1 + \vv{Y}_2 \\
  \del_0 &\dfn i \vv{Y}_3
\end{align*}
either as elements of $\C \su(2)$ or as complex left-invariant vector fields on $\SU(2)$. One checks by hand the following commutation relations
\begin{align*}
  [\vv{Y}_1, \vv{Y}_2] &= \vv{Y}_3 \\
  [\vv{Y}_2, \vv{Y}_3] &= \vv{Y}_1 \\
  [\vv{Y}_3, \vv{Y}_1] &= \vv{Y}_2
\end{align*}
(which characterize the Lie algebra $\su(2)$ up to isomorphism) and that the assignment
\begin{align}
  (A,B) &\longmapsto 2 \ \mathrm{tr}(A \cdot B^*) \label{eq:matrix_ip}
\end{align}
defines the unique Euclidean inner product on $\su(2)$ for which $\vv{Y}_1, \vv{Y}_2, \vv{Y}_3$ is an orthonormal basis: the commutation relations above ensure that~\eqref{eq:matrix_ip} is $\ad$-invariant. With this choice of metric, the associated Laplace-Beltrami operator is
\begin{align*}
  \Delta = -\left( \vv{Y}_1^2 + \vv{Y}_2^2 + \vv{Y}_3^2 \right) = \del_0^2 + \frac{1}{2} \left( \del_+ \del_- + \del_- \del_+ \right).
\end{align*}
For further reference we also take note of the commutation relations
\begin{align*}
  [\del_+, \del_-] &= 2 \del_0, \\
  [\del_+, \del_0] &= - \del_+, \\
  [\del_-, \del_0] &= \del_-.
\end{align*}
\subsection{Representations of $\SU(2)$}

For each $\ell \in \frac{1}{2} \Z_+$ let $V_\ell \sset \C[z_1, z_2]$ stand for the space of all homogeneous polynomials of degree $2\ell$. For $k \in \{0, \ldots, 2\ell\}$ we define the monomials
\begin{align*}
  p_{\ell k}(z) &\dfn z_1^k z_2^{2\ell - k}
\end{align*}
which form a basis for $V_\ell$, so $\dim_\C V_\ell = 2\ell + 1$. Such polynomials can be pushed forward by the map $\theta$ and integrated over $\SU(2)$ (which is compact and carries its Haar measure), hence the $L^2$ inner product on $\SU(2)$ induces an inner product on $V_\ell$. Explicitly:
\begin{align}
  (f,g) \in V_\ell \times V_\ell &\longmapsto \int_{\SU(2)} f(\theta^{-1}(x)) \ \overline{g(\theta^{-1}(x))} \ \dd \mu(x) = \left \langle f \circ \theta^{-1}, g \circ \theta^{-1} \right \rangle_{L^2(\SU(2))}. \label{eq:Vlip}
\end{align}
For $k \in \{ -\ell + j \st j = 0, 1, \ldots, 2\ell \}$ we shall write\footnote{Although both $\ell$ and $k$ are ``half-integers'' we have $\ell - k, \ell + k \in \{ 0, 1, \ldots, 2\ell \}$ -- thus integers.}
\begin{align}
  q_{\ell k}(z) \dfn \frac{z_1^{\ell - k} z_2^{\ell + k}}{\sqrt{(\ell - k)!(\ell + k)!}} = \frac{p_{\ell,\ell - k}(z)}{\sqrt{(\ell - k)!(\ell + k)!}} \label{eq:basis_qlk}
\end{align}
which form an orthonormal basis for $V_\ell$. We define $T_\ell: \SU(2) \rarr \GL(V_\ell)$ by\footnote{As usual, we are identifying a complex vector $z = (z_1, z_2) \in \C^2$ with the $2 \times 1$ matrix $\left( \begin{array}{c} z_1 \\ z_2 \end{array} \right)$.}
\begin{align*}
  \left( T_\ell (x) f \right) (z) &\dfn f(\transp{x} \cdot z).
\end{align*}
It turns out that $T_\ell$ is a group homomorphism i.e.~$(T_\ell, V_\ell)$ is a representation of $\SU(2)$, which can be shown to be irreducible and unitary w.r.t.~\eqref{eq:Vlip}. Moreover, if $(\xi, V_\xi)$ is any other irreducible unitary representation of $\SU(2)$ then $(\xi, V_\xi) \cong (T_\ell, V_\ell)$ for some $\ell \in \frac{1}{2} \Z_+$ i.e.
\begin{align*}
  \widehat{\SU(2)} &= \left\{ [T_\ell] \st \ell \in \frac{1}{2} \Z_+ \right\}.
\end{align*}

We shall write $t^\ell_{mn}$ (for $m,n \in \{ -\ell + j \st 0, 1, \ldots, 2\ell \}$) for the matrix elements of $(T_\ell, V_\ell)$ w.r.t.~the basis $\{ q_{\ell k} \}$ in~\eqref{eq:basis_qlk}. Therefore~\cite[eqn.~(2.5)]{araujo19} becomes
\begin{align}
  f &= \sum_{\ell \in \frac{1}{2} \Z_+} (2\ell + 1) \sum_{m,n = -\ell}^{\ell} \langle f, t^\ell_{mn} \rangle_{L^2(\SU(2))} \ t^\ell_{mn} \label{eq:fourier_su}.
\end{align}
Of fundamental importance to our ability to compute symbols of left-invariant operators effectively is that we know how our basic vector fields operate on the matrix elements~\cite[Proposition~11.9.2]{rt_psdoas}:
\begin{align*}
  \vv{Y}_1 t^\ell_{mn} &= \frac{\sqrt{(\ell - n)(\ell + n + 1)}}{-2i} t^\ell_{m,n + 1} + \frac{\sqrt{(\ell + n)(\ell - n + 1)}}{-2i} t^\ell_{m,n - 1}, \\
  \vv{Y}_2 t^\ell_{mn} &= \frac{\sqrt{(\ell - n)(\ell + n + 1)}}{2} t^\ell_{m,n + 1} - \frac{\sqrt{(\ell + n)(\ell - n + 1)}}{2} t^\ell_{m,n - 1}, \\
  \vv{Y}_3 t^\ell_{mn} &= -i n t^\ell_{mn}.
\end{align*}
Or, in another basis~\cite[Theorem~11.9.3]{rt_psdoas}:
\begin{align*}
  \del_+ t^\ell_{mn} &= - \sqrt{(\ell - n)(\ell + n + 1)} t^\ell_{m,n + 1}, \\
  \del_- t^\ell_{mn} &= - \sqrt{(\ell + n)(\ell - n + 1)} t^\ell_{m,n - 1}, \\
  \del_0 t^\ell_{mn} &= n t^\ell_{mn}.
\end{align*}
As for the Laplace-Beltrami operator we have\footnote{Notice that our Laplace-Beltrami differs from Ruzhansky and Turunen's (which they denote by $\mathcal{L}$) by a sign.}
\begin{align*}
  \Delta t^\ell_{mn} &= \ell(\ell + 1) t^\ell_{mn}.
\end{align*}
From this last fact we conclude that $\widehat{\SU(2)}$ is in one-to-one correspondence with $\sigma(\Delta)$:
\begin{align*}
  \sigma(\Delta) &= \left \{ \ell(\ell + 1) \st \ell \in \frac{1}{2} \Z_+ \right \}
\end{align*}
and
\begin{align*}
  E_{\ell(\ell + 1)} &= \mathcal{M}_{T_\ell}, \quad \forall \ell \in \frac{1}{2} \Z_+.
\end{align*}
In particular $\dim E_{\ell(\ell + 1)} = \dim \mathcal{M}_{T_\ell} = (2\ell + 1)^2$.

\section{An example: the standard CR structure on $\s^3$}

As a real hypersurface in $\C^2$, the sphere $\s^3$ carries a natural CR structure: for $p \in \s^3$ we denote by $\VV_p$ the space of complex tangent vectors of the form
\begin{align}
  a_1 \left. \frac{\del}{\del \bar{z}_1} \right|_p + a_2 \left. \frac{\del}{\del \bar{z}_2} \right|_p, \quad a_1, a_2 \in \C, \label{eq:antiholvec}
\end{align}
which happen to lie in $\C T_p \s^3$; this determines the fiber at $p$ of a vector subbundle $\VV \sset \C T \s^3$, which can be shown to be involutive and, actually, a \emph{CR structure} i.e.~$\VV_p \cap \bar{\VV}_p = 0$ for every $p \in \s^3$. Our map $\theta: \C^2 \rarr \gl(2, \C)$ previously defined pulls back a group structure from $\SU(2)$ to $\s^3$ which renders them isomorphic as Lie groups (via $\theta$):
\begin{align*}
  z \cdot w &\dfn \theta^{-1} \left( \theta(z) \cdot \theta(w) \right), \quad z, w \in \s^3.
\end{align*}
We will show that $\VV$ is left-invariant w.r.t.~this group structure; this is of course the same as saying that $\theta_* \VV$ is left-invariant on $\SU(2)$.

First of all, recall that since $\rho(z) \dfn |z|^2 - 1$ yields a smooth, globally defined defining function for $\s^3$ i.e.~$\s^3 = \rho^{-1}(0)$, we have that
\begin{align*}
  \C T_p \s^3 &= \{ v \in \C T_p \C^2 \st \dd \rho(v) = 0 \}, \quad \forall p \in \s^3.
\end{align*}
On the other hand, if $v \in \C T_p \C^2$ is of the form~\eqref{eq:antiholvec} then
\begin{align*}
  \dd \rho(v) = v(\rho) = a_1 \left. \frac{\del \rho}{\del \bar{z}_1} \right|_p + a_2 \left. \frac{\del \rho}{\del \bar{z}_2} \right|_p = a_1 p_1 + a_2 p_2 = \langle a, \bar{p} \rangle
\end{align*}
where $a \dfn (a_1, a_2)$. Therefore we have
\begin{align*}
  \VV_p &= \left\{ a_1 \left. \frac{\del}{\del \bar{z}_1} \right|_p + a_2 \left. \frac{\del}{\del \bar{z}_2} \right|_p \st \langle a, \bar{p} \rangle = 0 \right\}, \quad \forall p \in \s^3.
\end{align*}
If we take $a_1 \dfn -p_2$ and $a_2 \dfn p_1$ we have
\begin{align*}
  \VV_p &= \Span_\C \left\{ -p_2 \left. \frac{\del}{\del \bar{z}_1} \right|_p + p_1 \left. \frac{\del}{\del \bar{z}_2} \right|_p \right\}
\end{align*}
meaning that the vector field
\begin{align*}
  \LL &\dfn -z_2 \frac{\del}{\del \bar{z}_1} + z_1 \frac{\del}{\del \bar{z}_2}
\end{align*}
is a global section of $\VV$ on $\s^3$.

Next, we check that $\LL$ is a left-invariant vector field on $\s^3$. For $w = (w_1, w_2), z = (z_1, z_2) \in \C^2$ we have
\begin{align*}
  \theta(w_1,w_2) \cdot \theta(z_1,z_2) &= 
  \left(
    \begin{array}{c c}
      w_1 & -\bar{w}_2 \\
      w_2 & \bar{w}_1
    \end{array}
    \right) \cdot
  \left(
    \begin{array}{c c}
      z_1 & -\bar{z}_2 \\
      z_2 & \bar{z}_1
    \end{array}
    \right) \\
  &= \left(
    \begin{array}{c c}
      w_1 z_1 - \bar{w}_2 z_2 & -w_1 \bar{z}_2 - \bar{w}_2 \bar{z}_1 \\
      w_2 z_1 + \bar{w}_1 z_2 & -w_2 \bar{z}_2 + \bar{w}_1 \bar{z}_1
    \end{array}
    \right) \\
  &= \theta(w_1 z_1 - \bar{w}_2 z_2, w_2 z_1 + \bar{w}_1 z_2),
\end{align*}
hence
\begin{align*}
  w \cdot z &= (w_1 z_1 - \bar{w}_2 z_2, w_2 z_1 + \bar{w}_1 z_2), \quad \forall w,z \in \s^3.
\end{align*}
If $f \in \cinfty(\C^2)$ we have
\begin{align*}
  \LL (f \circ L_w)(z) &= \left( -z_2 \frac{\del}{\del \bar{z}_1} + z_1 \frac{\del}{\del \bar{z}_2} \right) f(w_1 z_1 - \bar{w}_2 z_2, w_2 z_1 + \bar{w}_1 z_2) \\
  &= -z_2 \left( \frac{\del f}{\del \bar{z}_1}(w \cdot z) \bar{w}_1 + \frac{\del f}{\del \bar{z}_2}(w \cdot z) \bar{w}_2 \right) + z_1 \left( \frac{\del f}{\del \bar{z}_1}(w \cdot z) (- w_2) + \frac{\del f}{\del \bar{z}_2}(w \cdot z) w_1 \right) \\
  &= -(w_2 z_1 + \bar{w}_1 z_2) \frac{\del f}{\del \bar{z}_1}(w \cdot z) + (w_1 z_1 - \bar{w}_2 z_2) \frac{\del f}{\del \bar{z}_2}(w \cdot z) \\
  &= (\LL f)(L_w z)
\end{align*}
for every $w$ and $f \in \cinfty(\C^2)$, thus proving that $\LL$ is left-invariant. We also note that
\begin{align*}
  \LL|_{(1,0)} = \left. \frac{\del}{\del \bar{z}_2} \right|_{(1,0)} = \frac{1}{2} \left( \left. \frac{\del}{\del x_2} \right|_{(1,0)} + i \left. \frac{\del}{\del y_2} \right|_{(1,0)} \right).
\end{align*}
Since by the identification $T_{(1,0)} \C^2 \cong \R^4 \cong \C^2$ we have
\begin{align*}
  \left. \frac{\del}{\del x_2} \right|_{(1,0)} &\cong (0,1) \\
  \left. \frac{\del}{\del y_2} \right|_{(1,0)} &\cong (0,i)
\end{align*}
hence in $\SU(2)$: 
\begin{align*}
  \theta_* \LL = \frac{1}{2} \theta(0,1) + \frac{i}{2} \theta(0,i) = \vv{Y}_2 + i \vv{Y}_1 = \del_-.
\end{align*}
\section{Structures of corank $2$ (i.e.~vector fields) on $\SU(2)$}

\subsection{Algebraic computations on $3$-dimensional Lie algebras: corank $2$ subalgebras}

For the moment, $\gr{g}$ will denote the Lie algebra of a general (compact) Lie group of dimension $3$. We will compute the algebraic expression of $\dd'$ (in non-trivial bidegrees) associated with Lie subalgebras $\gr{v} \sset \C \gr{g}$, and then employ these computations to study some examples of left-invariant structures on $\SU(2)$.

Suppose that $\gr{v}$ has dimension $n = 1$. We select a non-zero $\LL \in \gr{v}$, and $\MM_1, \MM_2 \in \C \gr{g}$ such that $\LL, \MM_1, \MM_2$ form a basis for $\C \gr{g}$: let $\tau, \zeta_1, \zeta_2 \in \C \gr{g}^*$ denote the corresponding dual basis. We write, for $j \in \{1, 2\}$,
\begin{align*}
  [\MM_1, \MM_2] &= a_1 \MM_1 + a_2 \MM_2 + a_3 \LL \\
  [\LL, \MM_j] &= b_{j,1} \MM_1 + b_{j,2} \MM_2 + b_{j,3} \LL
\end{align*}
and
\begin{align*}
  \dd \zeta_j &= d_{j,1} \zeta_1 \wedge \tau + d_{j,2} \zeta_2 \wedge \tau + d_{j,3} \zeta_1 \wedge \zeta_2. 
\end{align*}
Thus for $k \in \{1,2\}$ we have
\begin{align*}
  d_{j,k} = \dd \zeta_j (\MM_k, \LL) = \zeta_j ([\LL, \MM_k]) = \zeta_j \left( b_{k,1} \MM_1 + b_{k,2} \MM_2 + b_{k,3} \LL \right) = b_{k,j}
\end{align*}
while for $k = 3$ we have
\begin{align*}
  d_{j,3} = \dd \zeta_j (\MM_1, \MM_2) = \zeta_j ([\MM_2, \MM_1]) = \zeta_j \left( -a_1 \MM_1 - a_2 \MM_2 - a_3 \LL \right) = -a_j.
\end{align*}
We conclude that
\begin{align*}
  \dd \zeta_j &= b_{1,j} \zeta_1 \wedge \tau + b_{2,j} \zeta_2 \wedge \tau - a_j \zeta_1 \wedge \zeta_2 
\end{align*}
from which it follows easily that
\begin{align*}
  \dd (\zeta_1 \wedge \zeta_2) &= (- b_{1,1} - b_{2,2}) \zeta_1 \wedge \zeta_2 \wedge \tau.
\end{align*}

We have thus established generators for the spaces $\Lambda^{p,q}$ in the relevant (non-trivial) cases, and their corresponding differentials, as summarized below:
\begin{enumerate}
\item $p = 1$, $q = 0$:
  \begin{align*}
    \zeta_1 &\rightsquigarrow \dd' \zeta_1 = b_{1,1} \zeta_1 \wedge \tau + b_{2,1} \zeta_2 \wedge \tau \\
    \zeta_2 &\rightsquigarrow \dd' \zeta_2 = b_{1,2} \zeta_1 \wedge \tau + b_{2,2} \zeta_2 \wedge \tau
  \end{align*}
\item $p = 2$, $q = 0$:
  \begin{align*}
    \zeta_1 \wedge \zeta_2 &\rightsquigarrow \dd'(\zeta_1 \wedge \zeta_2) = (-b_{1,1} - b_{2,2}) \zeta_1 \wedge \zeta_2 \wedge \tau
  \end{align*}
\end{enumerate}
Notice how these expressions do \emph{not} depend on the coefficients that describe how $\MM_1$ and $\MM_2$ commute -- which we already knew, since $\dd'$ is an intrinsic property of the Lie algebra $\gr{v}$. Therefore:
\begin{enumerate}
\item $p = 1$, $q = 0$:
  \begin{align*}
    u = u_1 \zeta_1 + u_2 \zeta_2 &\Longrightarrow \dd' u = \left( -\LL u_1  + b_{1,1} u_1 + b_{1,2} u_2 \right) \zeta_1 \wedge \tau + \left( -\LL u_2  + b_{2,1} u_1 + b_{2,2} u_2 \right) \zeta_2 \wedge \tau
  \end{align*}
\item $p = 2$, $q = 0$:
  \begin{align*}
    u = v \ \zeta_1 \wedge \zeta_2 &\Longrightarrow \dd' u = \left(\LL v - (b_{1,1} + b_{2,2})v \right)  \zeta_1 \wedge \zeta_2 \wedge \tau
  \end{align*}
\end{enumerate}
\subsection{Example: analysis of the $\del_-$ operator}

Back to $\SU(2)$, we shall compute the action of $\del_-$ on each $\mathcal{M}_{T_\ell}$ where $\ell \in \frac{1}{2} \Z_+$. Let $\phi \in \mathcal{M}_{T_\ell}$ i.e.
\begin{align}
  \phi &= \sum_{m,n = -\ell}^\ell c_{mn} \sqrt{2\ell + 1} t^\ell_{mn} \label{eq:M_l}
\end{align}
for some constants $c_{mn} \in \C$. On passing, we recall that
\begin{align*}
  \| \phi \|_{L^2(\SU(2))} &= \left( \sum_{m,n = -\ell}^\ell |c_{mn}|^2 \right)^{\frac{1}{2}}.
\end{align*}
Suppose that $\del_- \phi = 0$, that is,
\begin{align*}
  \sum_{m,n = -\ell}^\ell c_{mn} \sqrt{2\ell + 1} \sqrt{(\ell + n)(\ell - n + 1)} t^\ell_{m, n - 1} &= 0.
\end{align*}
Then $c_{mn} = 0$ whenever $n \neq -\ell$, that is,
\begin{align*}
  \ker \del_- |_{\mathcal{M}_{T_\ell}} &= \span \left\{ t^\ell_{m,-\ell} \st m \in \frac{1}{2} \Z_+, \ -\ell \leq m \leq \ell \right \}.
\end{align*}
Therefore $\dim \ker \del_- |_{\mathcal{M}_{T_\ell}} = 2\ell + 1$; in particular $\del_- |_{\mathcal{M}_{T_0}} = 0$ (which we already knew since $\del_-$ is a vector field and $\mathcal{M}_{T_0}$ is the space of constants).
\begin{Prop} \label{prop:delminusclosed} The map $\del_-: \cinfty(\SU(2)) \rarr \cinfty(\SU(2))$ has closed range.
  \begin{proof} Suppose that $\phi$ in~\eqref{eq:M_l} is orthogonal to $\ker \del_- |_{\mathcal{M}_{T_\ell}}$ i.e.
    \begin{align*}
      \phi &= \sum_{m= -\ell}^\ell \sum_{n= -\ell+ 1}^\ell c_{mn} \sqrt{2\ell + 1} t^\ell_{mn}.
    \end{align*}
    Assume moreover that $ \| \phi \|_{L^2(\SU(2))} = 1$ (notice that this implies that $\ell \geq \frac{1}{2}$ since $\del_-|_{\mathcal{M}_{T_0}} = 0$). Then
    \begin{align*}
      \del_- \phi &= - \sum_{m= -\ell}^\ell \sum_{n= -\ell + 1}^\ell c_{mn} \sqrt{2\ell + 1} \sqrt{(\ell + n)(\ell - n + 1)} t^\ell_{m, n - 1} \\
    \end{align*}
    which implies that
    \begin{align*}
      \| \del_- \phi \|_{L^2(\SU(2))}^2 &= \sum_{m= -\ell}^\ell \sum_{n= -\ell+ 1}^\ell |c_{mn}|^2(\ell + n)(\ell - n + 1)
    \end{align*}
    Notice however that the function
    \begin{align*}
      x \in \R \longmapsto (\ell + x)(\ell - x + 1) = \ell^2 - x^2 + \ell + x
    \end{align*}
    has a unique critical point at $x = 1/2$, which is a global maximum: its minimum over $[-\ell + 1, \ell]$ must be attained at $x = -\ell + 1$ and/or $x = \ell$. In either case, its value is $2\ell$ i.e.
    \begin{align*}
      (\ell + x)(\ell - x + 1) &\geq 2\ell, \quad \forall x \in [-\ell + 1, \ell]
    \end{align*}
    hence
    \begin{align*}
      \| \del_- \phi \|_{L^2(\SU(2))}^2 \geq 2\ell \sum_{m= -\ell}^\ell \sum_{n= -\ell+ 1}^\ell |c_{mn}|^2 = 2\ell.
    \end{align*}
    We conclude that
    \begin{align*}
      \| \del_- \phi \|_{L^2(\SU(2))} &\geq \sqrt{2\ell}, \quad \forall \phi \in \mathcal{M}_{T_\ell} \cap \left(\ker \del_- \right)^\bot, \ \| \phi \|_{L^2(\SU(2))} = 1.
    \end{align*}
    On the other hand, there exists a constant $C > 0$ such that
    \begin{align*}
      \sqrt{2\ell} &\geq C (1 + \ell(\ell + 1))^\frac{1}{3}, \quad \forall \ell \geq \frac{1}{2},
    \end{align*}
    meaning that~\cite[eqn.~(6.1)]{araujo19} holds for $P = \del_-$ with $s = 1/3$, and $\omega = \omega_\infty$ as defined in~\cite[Sec~5.1]{araujo19}. By~\cite[Theorem~6.4(1)]{araujo19}, $\del_-$ is almost $\cinfty$ globally hypoelliptic, and this is equivalent, according to~\cite[Theorem~6.7(1)]{araujo19}, to $\del_-: \cinfty(\SU(2)) \rarr \cinfty(\SU(2))$ having a closed range.
  \end{proof}
\end{Prop}
In order to regard $\del_-$ as a left-invariant \emph{structure} on $\SU(2)$ -- and then to analyze the associated operator $\dd'$ acting on forms -- we must introduce some ``complementary'' vectors fields $\MM_1, \MM_2$. Of course one natural choice is to take
\begin{align*}
  \LL &\dfn \del_- \\
  \MM_1 &\dfn \del_+ \\
  \MM_2 &\dfn \del_0
\end{align*}
which yields the constants (in our previous notation)
\begin{align*}
  \begin{array}{|c c c|}
    \hline
    a_1 = -1, & a_2 = 0, & a_3 = 0 \\
    b_{1,1} = 0, & b_{1,2} = -2, & b_{1,3} = 0 \\
    b_{2,1} = 0, & b_{2,2} = 0, & b_{2,3} = 1 \\
    \hline
  \end{array}
\end{align*}
and therefore the operator $\dd'$ can be expressed, in the non-trivial bidegrees, as the matrices:
\begin{align*}
  \TR{\dd'_{(1,0)}}{(u_1, u_2)}{\cinfty(\SU(2))^2}{(- \del_- u_1 - 2 u_2, -\del_- u_2)}{\cinfty(\SU(2))^2},
\end{align*}
\begin{align*}
  \TR{\dd'_{(2,0)}}{u}{\cinfty(\SU(2))}{\del_- u}{\cinfty(\SU(2))}.
\end{align*}

Let us work on $\dd'_{(1,0)}$. For $\ell \in \frac{1}{2} \Z_+$, easy computations show that $(\phi_1, \phi_2) \in \mathcal{M}_{T_\ell}^2$ belongs to $\ker \dd'_{(1,0)}$ if and only if $\del_- \phi_2 = 0$ and $\del_- \phi_1 = -2 \phi_2$, that is:
\begin{itemize}
\item for $\ell = 0$, $\phi_2 = 0$ and $\phi_1 \in \mathcal{M}_{T_0}$;
\item for $\ell \geq \frac{1}{2}$ we have
  \begin{align*}
    \phi_2 &= \sqrt{2\ell + 1} \sum_{m = -\ell}^\ell c^2_{m, -\ell} t^\ell_{m, -\ell} \\
    \phi_1 &= \sqrt{2\ell + 1} \sum_{m = -\ell}^\ell \left( c^1_{m, -\ell}  t^\ell_{m, -\ell} + \frac{2}{\sqrt{2\ell}} c^2_{m, -\ell} t^\ell_{m, -\ell + 1} \right)
  \end{align*}
  where $c^1_{m, -\ell}, c^2_{m, -\ell} \in \C$.
\end{itemize}
In particular:
\begin{align*}
  \dim \ker \dd'_{(1,0)}|_{\mathcal{M}_{T_\ell}^2} &=
  \begin{cases}
    1, & \text{if $\ell = 0$}; \\
    2(2 \ell + 1), & \text{if $\ell \geq 1/2$}.
  \end{cases}
\end{align*}
\begin{Cor} The map $\dd'_{(1,0)}: \cinfty(\SU(2))^2 \rarr \cinfty(\SU(2))^2$ has closed range.
  \begin{proof} Let $ \{ (u_{1,\nu}, u_{2, \nu})\}_{\nu \in \N} \sset \cinfty(\SU(2))^2$ be a sequence and $(v_1, v_2) \in \cinfty(\SU(2))^2$ be such that, as $\nu \to \infty$,
    \begin{align*}
      \dd'_{(1,0)} (u_{1,\nu}, u_{2, \nu}) &\longrightarrow (v_1, v_2) \ \text{in $\cinfty(\SU(2))^2$},
    \end{align*}
    that is:
    \begin{align*}
      - \del_- u_{1, \nu} - 2 u_{2, \nu} &\longrightarrow v_1 \ \text{in $\cinfty(\SU(2))$}, \\
      -\del_- u_{2, \nu} &\longrightarrow v_2 \ \text{in $\cinfty(\SU(2))$}.
    \end{align*}
    Since $\del_-: \cinfty(\SU(2)) \rarr \cinfty(\SU(2))$ has closed range (Proposition~\ref{prop:delminusclosed}) there exists $u_2 \in \cinfty(\SU(2))$ such that $u_{2, \nu} \to u_2$ in $\cinfty(\SU(2))$. From this we conclude that $- \del_- u_{1, \nu} \to v_1 + 2 u_2$ in $\cinfty(\SU(2))$, which in turn implies that there exists $u_1 \in \cinfty(\SU(2))$ such that $u_{1, \nu} \to u_1$ in $\cinfty(\SU(2))$. This finishes the proof.
  \end{proof}
\end{Cor}
\section{Structures of corank $1$ on $\SU(2)$}

\subsection{Algebraic computations on $3$-dimensional Lie algebras: corank $1$ subalgebras}

Back to our general algebraic computations, let $\gr{v} \sset \C \gr{g}$ be a Lie subalgebra of complex dimension $n = 2$, where $\gr{g}$ is the Lie algebra of some $3$-dimensional compact Lie group. It should be pointed out that, due to dimensional reasons, such a subalgebra is either elliptic (i.e.~$\gr{v} + \bar{\gr{v}} = \C \gr{g})$ or essentially real (i.e.~$\gr{v} = \bar{\gr{v}}$). We select $\LL_1, \LL_2$ a basis of $\gr{v}$, and $\MM \in \C \gr{g}$ another vector field such that $\LL_1, \LL_2, \MM$ form a basis for $\C \gr{g}$: let $\tau_1, \tau_2, \zeta \in \C \gr{g}^*$ denote the corresponding dual basis. First, we write
\begin{align*}
  [\LL_1, \LL_2] &= a_1 \LL_1 + a_2 \LL_2 \\
  [\LL_j, \MM] &= b_{j,1} \LL_1 + b_{j,2} \LL_2 + b_{j,3} \MM
\end{align*}
hence
\begin{align*}
  \dd \zeta &= c_1 \zeta \wedge \tau_1 + c_2 \zeta \wedge \tau_2 \\
  \dd \tau_j &= d_{j,1} \zeta \wedge \tau_1 + d_{j,2} \zeta \wedge \tau_2 + d_{j,3} \tau_1 \wedge \tau_2 
\end{align*}
where
\begin{align*}
  c_j = \dd \zeta (\MM, \LL_j) = \zeta ([\LL_j, \MM]) = \zeta \left( b_{j,1} \LL_1 + b_{j,2} \LL_2 + b_{j,3} \MM \right) = b_{j,3}
\end{align*}
and for $k \in \{1,2\}$
\begin{align*}
  d_{j,k} = \dd \tau_j (\MM, \LL_k) = \tau_j ([\LL_k, \MM]) = \tau_j \left( b_{k,1} \LL_1 + b_{k,2} \LL_2 + b_{k,3} \MM \right) = b_{k,j}
\end{align*}
while for $k = 3$
\begin{align*}
  d_{j,3} = \dd \tau_j (\LL_1, \LL_2) = \tau_j ([\LL_2, \LL_1]) = \tau_j \left( -a_1 \LL_1 - a_2 \LL_2 \right) = -a_j.
\end{align*}
We conclude that
\begin{align*}
  \dd \zeta &= b_{1,3} \zeta \wedge \tau_1 + b_{2,3} \zeta \wedge \tau_2 \\
  \dd \tau_j &= b_{1,j} \zeta \wedge \tau_1 + b_{2,j} \zeta \wedge \tau_2 - a_j \tau_1 \wedge \tau_2 
\end{align*}
from which it follows easily that
\begin{align*}
  \dd (\zeta \wedge \tau_1) &= (a_1 - b_{2,3}) \zeta \wedge \tau_1 \wedge \tau_2, \\
  \dd (\zeta \wedge \tau_2) &= (a_2 + b_{1,3}) \zeta \wedge \tau_1 \wedge \tau_2.
\end{align*}

We have thus established generators for the spaces $\Lambda^{p,q}$ in the relevant (non-trivial) cases, and their corresponding differentials, as summarized below:
\begin{enumerate}
\item $p = 0$, $q = 1$:
  \begin{align*}
    \tau_1 &\rightsquigarrow \dd' \tau_1 = -a_1 \tau_1 \wedge \tau_2 \\
    \tau_2 &\rightsquigarrow \dd' \tau_2 = -a_2 \tau_1 \wedge \tau_2
  \end{align*}
\item $p = 1$, $q = 0$:
  \begin{align*}
    \zeta &\rightsquigarrow \dd' \zeta = b_{1,3} \zeta \wedge \tau_1 + b_{2,3} \zeta \wedge \tau_2 \\
  \end{align*}
\item $p = 1$, $q = 1$:
  \begin{align*}
    \zeta \wedge \tau_1 &\rightsquigarrow \dd' (\zeta \wedge \tau_1) = (a_1 - b_{2,3}) \zeta \wedge \tau_1 \wedge \tau_2 \\
    \zeta \wedge \tau_2 &\rightsquigarrow \dd' (\zeta \wedge \tau_2) = (a_2 + b_{1,3}) \zeta \wedge \tau_1 \wedge \tau_2
  \end{align*}
\end{enumerate}

Therefore:
\begin{enumerate}
\item $p = 0$, $q = 1$:
  \begin{align*}
    u = u_1 \tau_1 + u_2 \tau_2 &\Longrightarrow \dd' u = \left( \LL_1 u_2 - \LL_2 u_1 - a_1 u_1 - a_2 u_2 \right) \tau_1 \wedge \tau_2
  \end{align*}
\item $p = 1$, $q = 0$:
  \begin{align*}
    u = v \zeta &\Longrightarrow \dd' u = (- \LL_1 v + b_{1,3} v) \zeta \wedge \tau_1 + (- \LL_2 v + b_{2,3} v) \zeta \wedge \tau_2
  \end{align*}
\item $p = 1$, $q = 1$:
  \begin{align*}
    u = u_1 \ \zeta \wedge \tau_1 + u_2 \ \zeta \wedge \tau_2 &\Longrightarrow \dd' u = \left( \LL_2 u_1 - \LL_1 u_2 + (a_1 - b_{2,3})u_1 + (a_2 + b_{1,3})u_2 \right) \zeta \wedge \tau_1 \wedge \tau_2 
  \end{align*}
\end{enumerate}

\subsection{Example: a corank $1$ structure on $\SU(2)$}

In $\SU(2)$, let $\gr{v} \dfn \Span \{ \del_-, \del_0 \}$ and take
\begin{align*}
  \LL_1 &\dfn \del_- \\
  \LL_2 &\dfn \del_0 \\
  \MM &\dfn \del_+
\end{align*}
from which we have, in the notation of the previous section,
\begin{align*}
  \begin{array}{|c c c|}
    \hline
    \multicolumn{3}{|c|}{
      \begin{array}{c c}
        a_1 = 1, & a_2 = 0
      \end{array}
    } \\
    b_{1,1} = 0, & b_{1,2} = -2, & b_{1,3} = 0 \\
    b_{2,1} = 0, & b_{2,2} = 0, & b_{2,3} = -1 \\
    \hline
  \end{array}
\end{align*}
which we can use to express $\dd'$ in any non-trivial bidegree:
\begin{align*}
  \TR{\dd'_{(0,1)}}{(u_1, u_2)}{\cinfty(\SU(2))^2}{\del_- u_2 - \del_0 u_1 - u_1}{\cinfty(\SU(2))},
\end{align*}
\begin{align*}
  \TR{\dd'_{(1,0)}}{u}{\cinfty(\SU(2))}{(-\del_- u, -\del_0 u - u)}{\cinfty(\SU(2))^2},
\end{align*}
\begin{align*}
  \TR{\dd'_{(1,1)}}{(u_1, u_2)}{\cinfty(\SU(2))^2}{\del_0 u_1 - \del_- u_2 + 2 u_1}{\cinfty(\SU(2))}.
\end{align*}
We remark that $\gr{v}$ is elliptic since
\begin{align*}
  \bar{\del}_- = \overline{i \vv{Y}_1 + \vv{Y}_2} = -i \vv{Y}_1 + \vv{Y}_2 = - \del_+.
\end{align*}
We denote by $\VV \sset \C T \SU(2)$ the underlying left-invariant involutive structure.

\begin{Thm} If $\gr{v} \dfn \Span \{ \del_-, \del_0 \}$ then $H_{\VV}^{1,1} \left( \SU(2); \cinfty(\SU(2) \right) = 0$.
  \begin{proof} We claim that for each $\ell \in \frac{1}{2} \Z_+$ we have
    \begin{align*}
      \ran \left\{ \dd'_{(1,0)}: \mathcal{M}_{T_\ell} \longrightarrow \mathcal{M}_{T_\ell}^2 \right\} &= \ker \left\{ \dd'_{(1,1)}: \mathcal{M}_{T_\ell}^2 \longrightarrow \mathcal{M}_{T_\ell} \right\}.
    \end{align*}
    Indeed, notice that:
    \begin{itemize}
    \item For each $\ell \in \frac{1}{2} \Z_+$ the map $\dd'_{(1,0)}: \mathcal{M}_{T_\ell} \rarr \mathcal{M}_{T_\ell}^2$ is injective. Indeed for any $\phi \in \mathcal{M}_{T_\ell}$, which we again represent as in~\eqref{eq:M_l}, such that $\dd'_{(1,0)} \phi = (0,0)$, we have, in particular, that
      \begin{align*}
        0 = \del_0 \phi + \phi = \sum_{m,n = -\ell}^\ell c_{mn} \sqrt{2\ell + 1} (n + 1) t^\ell_{mn}
      \end{align*}
      which implies that the constants $c_{mn}$ are all zero: either $\ell = 0$ (hence $n = 0$), or $n + 1$ is \emph{never} zero (because in that case $n$ is a half-integer). Thus $\phi = 0$.
    \item For each $\ell \in \frac{1}{2} \Z_+$ the map $\dd'_{(1,1)}: \mathcal{M}_{T_\ell}^2 \rarr \mathcal{M}_{T_\ell}$ is surjective. Indeed, let $m,n \in [-\ell, \ell]$ be half-integers (or zero if $\ell = 0$) and take $\phi \dfn (n + 2)^{-1} t^\ell_{mn} \in \mathcal{M}_{T_\ell}$: clearly $\dd'_{(1,1)} (\phi, 0) = t^\ell_{mn}$.
    \item By the Rank-Nullity Theorem (we are working with finite-dimensional vector spaces!) the sequence
      \begin{align*}
        \xymatrix{
          \mathcal{M}_{T_\ell} \ar[r]^{\dd'_{(1,0)}} & \mathcal{M}_{T_\ell}^2 \ar[r]^{\dd'_{(1,1)}} &\mathcal{M}_{T_\ell}
        }
      \end{align*}
      is exact, thus proving our claim.
    \end{itemize}

    Since $\del_-: \cinfty(\SU(2)) \rarr \cinfty(\SU(2))$ has a closed range, it follows immediately from the expression of $\dd'_{(1,0)}: \cinfty(\SU(2)) \rarr \cinfty(\SU(2))^2$ that this map also has a closed range. This, together with our previous arguments, implies the existence of an isomorphism~\cite[Theorem~7.2]{araujo19}
    \begin{align*}
      H_{\VV}^{1,1} \left( \SU(2); \cinfty(\SU(2) \right) \cong \bigoplus_{\ell \in \frac{1}{2} \Z_+} \frac{ \ker \left\{ \dd'_{(1,1)}: \mathcal{M}_{T_\ell}^2 \longrightarrow \mathcal{M}_{T_\ell} \right\} }{ \ran \left\{ \dd'_{(1,0)}: \mathcal{M}_{T_\ell} \longrightarrow \mathcal{M}_{T_\ell}^2 \right\} } = 0.
    \end{align*}
  \end{proof}
\end{Thm}

\begin{Thm} If $\gr{v} \dfn \Span \{ \del_-, \del_0 \}$ then $\dim H_{\VV}^{0,1} \left( \SU(2); \cinfty(\SU(2) \right) = 1$.
  \begin{proof} First of all, notice that the map $\dd'_{(0,1)}: \mathcal{M}_{T_\ell}^2 \rarr \mathcal{M}_{T_\ell}$ is surjective for each $\ell \in \frac{1}{2} \Z_+$. Indeed, let $m,n \in [-\ell, \ell]$ be half-integers (or zero if $\ell = 0$) and take $\phi \dfn -(n + 1)^{-1} t^\ell_{mn} \in \mathcal{M}_{T_\ell}$: clearly $\dd'_{(0,1)} (\phi, 0) = t^\ell_{mn}$.

    Next, if $\ell \in \frac{1}{2} \Z_+$ is non-zero then the map $\dd'_{(0,0)}: \mathcal{M}_{T_\ell} \rarr \mathcal{M}_{T_\ell}^2$ is injective. Indeed, recall that $\dd'_{(0,0)} u = (\del_- u, \del_0 u)$: since $\del_0|_{\mathcal{M}_{T_\ell}}$ is clearly injective when $\ell \neq 0$ we have injectivity of $\dd'_{(0,0)}|_{\mathcal{M}_{T_\ell}}$ in that case. On the other hand, both $\del_-$ and $\del_0$ act as the zero map on $\mathcal{M}_{T_0}$ hence $\dd'_{(0,0)}|_{\mathcal{M}_{T_0}} = 0$.

    Therefore:
    \begin{align*}
      \bigoplus_{\ell \in \frac{1}{2} \Z_+} \frac{ \ker \left\{ \dd'_{(0,1)}: \mathcal{M}_{T_\ell}^2 \longrightarrow \mathcal{M}_{T_\ell} \right\} }{ \ran \left\{ \dd'_{(0,0)}: \mathcal{M}_{T_\ell} \longrightarrow \mathcal{M}_{T_\ell}^2 \right\} } &= \frac{ \ker \left\{ \dd'_{(0,1)}: \mathcal{M}_{T_0}^2 \longrightarrow \mathcal{M}_{T_0} \right\} }{ \ran \left\{ \dd'_{(0,0)}: \mathcal{M}_{T_0} \longrightarrow \mathcal{M}_{T_0}^2 \right\} }
      \end{align*}
    is one-dimensional by (several applications of) the Rank-Nullity Theorem.

    But, as in the proof of the previous theorem, $\dd'_{(0,0)}: \cinfty(\SU(2)) \rarr \cinfty(\SU(2))^2$ clearly has a closed range, hence
    \begin{align*}
      H_{\VV}^{0,1} \left( \SU(2); \cinfty(\SU(2) \right) &\cong \frac{ \ker \left\{ \dd'_{(0,1)}: \mathcal{M}_{T_0}^2 \longrightarrow \mathcal{M}_{T_0} \right\} }{ \ran \left\{ \dd'_{(0,0)}: \mathcal{M}_{T_0} \longrightarrow \mathcal{M}_{T_0}^2 \right\} }.
    \end{align*}
  \end{proof}
\end{Thm}


\def\cprime{$'$}

\end{document}